\documentclass[twoside,12pt]{amsart}
\usepackage[a4paper,scale=.795,centering]{geometry}  
\usepackage{mathrsfs,amssymb,mathabx, amsfonts,  amsthm, amscd, amsmath}
\usepackage[mathlines]{lineno}
\usepackage[dvipsnames]{xcolor}

\usepackage{graphicx}

\newtheorem{theorem}{Theorem}[section]

\newtheorem{lemma}[theorem]{Lemma}
\newtheorem{proposition}[theorem]{Proposition}
\newtheorem{remark}[theorem]{Remark}

\newtheorem{model}[theorem]{Model}

\newcommand{\VEC}[1]{\mathbf{#1}}
\newcommand{\mat}[1]{\mathbf{#1}}
\newcommand{\abs}[1]{\hspace{.25mm}\left|#1\right|\hspace{.25mm}}
\newcommand{\z}{\mathbb{Z}}
\newcommand{\imf}[2]{\ensuremath{#1\!\paren{#2}}}
\newcommand{\sgn}{\mathrm{sgn}}
\newcommand{\diag}{\mathrm{diag}}
\newcommand{\ud}{\mathrm{d}}
\newcommand{\paren}[1]{\left( #1\right) }
\newcommand{\na}{\mathbb{N}}

\newcommand{\re}{\mathbb{R}}

\newcommand{\set}[1]{\left\{ #1\right\} }

\newcommand{\indi}[1]{\si_{#1}}
\newcommand{\si}{{\mathbf{1}}}
\newcommand{\bra}[1]{\left[ #1\right] }
\newcommand{\se}{\mathbb{E}}

\title{On multitype Branching Processes with Interaction}
\author{Mar\'ia Clara Fittipaldi, Sandra Palau}
\address[MCF]{Facultad de Ciencias,  
	Universidad Nacional Aut\'onoma de M\'exico,   M\'exico. email: mcfittipaldi@ciencias.unam.mx}
\address[SP]{Instituto de Investigaciones en Matem\'aticas Aplicadas y en Sistemas,  	Universidad Nacional Aut\'onoma de M\'exico,   M\'exico. email: sandra@sigma.iimas.unam.mx }

\begin{document}
\begin{abstract}
Motivated by the stochastic Lotka-Volterra model, we introduce discrete-state interacting multitype branching processes. We show that they can be obtained as the sum of a multidimensional random walk with a Lamperti-type change proportional to the population size; and a multidimensional Poisson process with a time-change proportional to the pairwise interactions. We define the analogous continuous-state process as the unique strong solution of a multidimensional SDE. We prove that the scaling limits of the discrete-state process correspond to its continuous counterpart. In addition, we show that the continuous-state model can be constructed as a generalized Lamperti-type transformation of multidimensional Lévy processes.\\

\textbf{MSC2020 subject classifications:} 
 60F17
; 60G51
; 60H20
; 60J25
; 60J80
; 92D25
.

\end{abstract}
\maketitle

	\section{Introduction}
	
	In this article, we present a continuous-time and discrete-state interacting multitype  branching process, 
	where the interaction is intratype or intertype and corresponds to competition or cooperation. 
	We prove that it can  be obtained as the sum of  two  time-changed processes;  a multidimensional random walk with a Lamperti-type change proportional to the population size and  a multidimensional Poisson process with a time change proportional to the pairwise interactions.  The Markovian population dynamics associated with this process is the following
	
	\begin{model}\label{Mo: DIMBP}  \textbf{\textsc{Discrete-space interacting multitype  branching process (DIMBP)}} Consider a po\-pu\-lation with $d$ different types of individuals. 
		When there are $u_1,\ldots, u_d$ individuals of types $1,\ldots, d$, 
		each individual  of type $i\in \{1,\dots, d\}$ dies at rate $\lambda_i$ and leaves behind  $v_1,\ldots, v_d$ offspring of types $1,\ldots, d$ with probability $\mu_i(v_1,\ldots, v_d)$, independently of others. 
		Additionally, each individual of type $i$ chooses at rate $\abs{c^{i,j}}$, an individual of type $j$. If $c^{i,j}<0$, the chosen individual is killed. Otherwise an independent replica of the chosen individual is incorporated into the population. 
	\end{model}
	
	In the discrete-state setting, there are some related models. For example, in discrete time, we can interpret the dynamic of our process as a particular case of a population size-dependent multitype branching process studied in Gonzalez et. al. \cite{GMM2004,GMM2005}. For continuous-time processes, Champagnat and Villemonais \cite{ChV2015} presented a model where an individual generates one child of the same type with birth rate dependent on the whole population, and there is a non-linear interspecific and intraspecific competition.   
	
	\smallskip  
	
	In this manuscript, we define the continuous-state analogous model as the unique strong solution of a stochastic differential equation. We call it continuous-state interacting multitype branching process (CIMBP). 
	Here, the continuum mass of each type reproduces within its own population type and allows seeding of mass into other population types. In addition, there are both intra and intertype pairwise interactions  proportional to the product of their type-population masses. 
	
	CIMPs generalize the Lotka-Volterra model studied by  Cattiaux and M\'el\'eard, \cite{MR2606515} which consists in two independent continuous-state branching processes with continuous trajectories that have both pairwise interactions. There are other related models, for example \cite{MR2466448, HN2018}. 
	
	Our main theorem is the following: every CIMBP can be constructed as the scaling limit of a sequence of renormalized DIMBPs. 
	This result is obtained via a sequence of piecewise multitype branching processes approximation on a time-population grid. 
	As a consequence of the scaling limit representation, a CIMBP can be seen as a generalized Lamperti-type transformation of multidimensional L\'evy processes.
	
	The remainder of the paper is structured as follows. In Section \ref{discreto}, we study the discrete-state model and its Lamperti-type transformation. In Section \ref{continuo}, we define the continuous-state model in terms of a stochastic differential equation and we present a time-change equation associated with it. Our main theorem is presented at Subsection \ref{limit} and it is proved in Section \ref{se: Grid}. For simplicity of the exposition, the proofs of auxiliary results are presented in Section \ref{se: proofs}.
	
	\section{The discrete space model}\label{discreto}
	Consider a DIMBP $(\VEC{Z}_t, t\geq 0)$ associated with Model \ref{Mo: DIMBP}, i.e. a $\z^{d}_+$-valued Markovian process with the following dynamics
	\begin{itemize}
		\item the transition rate from $\VEC{u}$ to $\VEC{u}-\VEC{e_i}+\VEC{v}$ following a reproduction event of the $i$-th type is equal to $\lambda_i u_i\imf{\mu_i}{\VEC{v}}$; 
		\item the transition rate from $\VEC{u}$ to $\VEC{u}+\imf{\sgn}{c^{i,j}}\VEC{e_j}$ following an interaction event from type $i$ to type $j$ is equals to $\abs{c^{i,j}}u_iu_j$. 
	\end{itemize} 
	
	To be clear, for each $i\in D:=\{1,\dots,d\}$, $\mu_i$ is a distribution on $\na^d$, $\lambda_i>0$ and $c^{i,j}\in \re$.
	Without loss of generality, we now assume that $\mu_i(\VEC{e_i})=0.$
	
	We let $\mat{C}$ be the matrix $\left(c^{i,j},i,j \in D\right)$.  
	Note that if $\mat{C}$ is the zero matrix, we obtain a multitype continuous-time Galton-Watson process. 
	If $\mat{C}$ has negative entries, our model is a multitype extension of the logistic branching process introduced in \cite{MR2134113}. The multitype extension includes both intratype and intertype pairwise competition (leading to quadratic intratype competition). We obtain the discrete version of the predator-prey (or Lotka-Volterra) branching processes studied in \cite{MR2606515} when $d=2$,   $\mat{C}$ has zero diagonal, the off-diagonal entries have different signs, and there is no intratype reproduction.  
	
	\smallskip
	Let us turn to a construction of the above DIMBP in terms of random walks and Poisson processes. Our construction features multiparameter time-changes  as the ones introduced in \cite{MR577310} and used in the branching process setting in \cite{MR3444314}. 
	
	Consider $d$ independent random walks $\VEC{X^1},\ldots, \VEC{X^d}$ with values in $\z^d$,  where the jump rate of $\VEC{X^i}$ is $\lambda_i\tilde \mu_i$, with $\tilde \mu_i(\VEC{v})=\mu_i(\VEC{v}+\VEC{e_i})$. 
	Let $X^{i,j}$ be the $j$-entry of $\VEC{X^i}$. Note that $X^{i,i}$ is downwards skip-free, i.e. its jumps belong to $\set{-1,1,2, \ldots}$. In addition, $X^{i,j}$ is non-decreasing for $j\neq i$. 
	Let $N^{i,j}$  be unit rate Poisson processes. Suppose that all the processes are independent.
	Consider the process $\VEC{Z}=\left\{(Z_t^1,\ldots, Z_t^d), t\geq 0\right\}$ defined as 
	\begin{equation}
	\label{eq: DIMBP TC}
	Z^{j}_t=z^j +\sum_{i=1}^d X^{i,j}_{\int_0^tZ^{i}_r\ud r}+\sum_{i=1}^d  sgn(c^{i,j})N^{i,j}_{|c^{i,j}|\int_0^t Z^{i}_r Z^{j}_r \ud r},
	\qquad j \in D,
	\end{equation}
	with starting value $\VEC{z}=(z^i,i\in D) \in \z^d_+$.
	
	Observe that  the vector-valued  equation has to be solved simultaneously for all $j\in D$. The existence and uniqueness of a solution follows from the piecewise constant paths of  $X^{i,j}$ and $N^{i,j}$,  $i,j\in D$, which implies the same for every $Z^{j}$.
	
	\begin{remark}
		When $Z^j_t=0$, the values $X^{j,j}_{\int_0^tZ^{j}_r\ud r}$ and $N^{i,j}_{|c^{i,j}|\int_0^t Z^{i}_r Z^{j}_r \ud r}$ become constant and since $X^{i,j}$, $i\neq j$ are non-decreasing, then $Z^{j}_{t+\cdot}$ is non-negative. Moreover,  the vector zero is an absorbing point for $\VEC{Z}$. 
	\end{remark}
	
	\begin{proposition}
		The stochastic process $\VEC{Z}$ given in \eqref{eq: DIMBP TC} is a DIMBP. 
	\end{proposition}
	\begin{proof}
		The reader can note that any solution of equation \eqref{eq: DIMBP TC} is non-negative with piecewise constant paths.  By analyzing the constancy intervals of $\VEC{Z}$, one can see that $\VEC{Z}$ is Markovian and has the correct jump rates.
		In particular, the transition from $\VEC{u}$ to $\VEC{u}-\VEC{e_i}+\VEC{v}$ following a reproduction event of the $i$-th type is governed by the jumps of vector $\VEC{X^i}$; while the transition from $\VEC{u}$ to $\VEC{u}+\imf{\sgn}{c^{i,j}}\VEC{e_j}$ following an interaction event from type $i$ to type $j$ is governed by $N^{i,j}$.
	\end{proof}
	
	\section{Continuous-state model}\label{continuo}
	\smallskip
	
	For every $i\in D$, let $W^{(i)}_t$ be a standard Brownian motion and let $\mathcal{N}^{(i)}(\ud s,\ud \VEC{r},\ud u)$ be a Poisson random measure on $\re_+^{d+2}$ with intensity measure $\ud s m^{(i)}(\ud \VEC{r})\ud u$, where $m^{(i)}$ is a Borel measure on $\re_+^d\setminus \{0\}$ satisfying 
	\begin{equation}\label{eq: integral conditions}
	\int_{\re_+^d\setminus \{0\}}\left[||\VEC{r}|| \wedge ||\VEC{r}||^2 + \sum\limits_{j=1, j\neq i}^d r_j\right]m^{(i)}(\ud \VEC{r}) <\infty. 
	\end{equation}
	Denote by $\widetilde{\mathcal{N}}^{(i)}$ the compensated measure of $\mathcal{N}^{(i)}$. We assume that all the processes are
	independent.  Let $\boldsymbol{\sigma}=(\sigma^i, i \in D) \in \re_+^d$,  $\mat{C}=(c^{ij}, i,j \in D)\in \re^{d\times d}$, and $\mat{B}=(b^{ij}, i,j \in D)\in \re_{(+)}^{d\times d}$, where $\re_{(+)}^{d\times d}$ is the set of essentially non-negative $d\times d$ matrices, i.e. matrices with nonnegative off-diagonal entries.
	
	Let us consider the $d$-dimensional stochastic differential equation (SDE for short) 
	\begin{equation}
	\begin{split}
	\label{continuous MLB}
	Y^j_t&=y^j+ \sum\limits_{i=1}^d\int_0^t c^{ij}Y^i_sY^j_s\ud s +\sum\limits_{i=1}^d\int_0^t b^{ij}Y^i_s\ud s + \int_0^t \sqrt{2\sigma^jY^j_s}\ud W^{(j)}_s\\ 
	& +\sum\limits_{i=1}^d \int_0^t\int_{\re_+^d}\int_0^{\infty}r_j\indi{\{u\leq Y^i_{s-}\}} \widetilde{\mathcal{N}}^{(i)}(\ud s,\ud \VEC{r},\ud u)
	\qquad j \in D, t\geq 0,
	\end{split}
	\end{equation}
	with starting value $\VEC{y}=(y^i, i \in D) \in \re_+^d$. 
	
	\begin{remark}
		If $\mat{C}=\mat{0}$, there exists a unique strong solution which corresponds to a multitype continuous-state branching process (MCBP for short) studied by \cite{MR3340375}. If $\mat{C}\neq \mat{0}$, the associated non-linear terms can be interpreted as pairwise interactions.
	\end{remark}
	
	We define a {\bf continuous-state interacting multitype branching process (CIMBP)} as the vector-valued Markov process $\VEC{Y}=\{(Y^{1}_t,\ldots, Y^d_t), t\geq 0\}$ that is the unique strong solution to \eqref{continuous MLB}.
	The following theorem guarantees the existence of such solution. It is an extension of the papers by Barczy, Li and Pap \cite{MR3340375} and Ma \cite{MR3208120}. Its proof can be found in Section \ref{se: proofs}.
	
	\begin{theorem}\label{thm: continuous MLB}
		There exists a unique non-negative strong solution to \eqref{continuous MLB}. Moreover, the associated process $\VEC{Y}=\left\{(Y^1_t,\ldots, Y^d_t), t \in \re_+\right\}$ has infinitesimal generator  
		\begin{eqnarray*}
			\mathcal{A}\VEC{F}(\VEC{x})&=&\langle \diag(\VEC{C}\VEC{x}\VEC{x}^{T}) + \VEC{B}\VEC{x},\VEC{F}'(\VEC{x})\rangle +\sum\limits_{i=1}^d \sigma^ix_iF''_{ii}(\VEC{x})\\
			&&+\sum\limits_{i=1}^d x_i\int_{\re_+^d} \left[\VEC{F}(\VEC{x}+\VEC{r})-\VEC{F}(\VEC{x})-\langle\VEC{r},\VEC{F}'(\VEC{r})\rangle\right]m^{(i)}(\ud \VEC{r}),
		\end{eqnarray*}
		for $\VEC{F}\in C^{2}_{c}(\re^d_+,\re)$ and $\VEC{x} \in \re^d_+$,  where $F'_i$ and $F''_{ii}$, $i \in D$, denote the first  and second order partial derivatives of $\VEC{F}$ with respect to its $i$-th variable, respectively, and  $\VEC{F}'(\VEC{x})=(F'_1(\VEC{x}), \dots, F'_d(\VEC{x}))^T.$ 
	\end{theorem}

	\smallskip
	We recall that, by equation \eqref{eq: DIMBP TC}, a discrete-space interacting multitype branching process can be defined in terms of random walks and Poisson processes.  
	The scaling limits of random walks are known to be  L\'evy processes and, by the strong law of large numbers, the scaling limits of Poisson processes are deterministic L\'evy processes.
	Hence, it is natural to think that the scaling limits of DIMBPs can be constructed in terms of L\'evy processes.

	Let $\VEC{X^1},\ldots, \VEC{X^d}$ be independent $\re^d$-valued L\'evy processes. We suppose that $X^{i,i}$ has no negative jumps and that $X^{i,j}$ is a subordinator for $i\neq j \in  D$. Let $\VEC{Z}=\left\{(Z_t^1,\ldots Z_t^d), t\geq 0\right\} $ be a solution to the equation
	\begin{equation}
	\label{eq: CIMBP TC}
	Z^j_t
	=z^j+\sum_{i=1}^d \bra{	X^{i,j}_{\int_0^t Z^i_s\, \ud s}+c^{i,j} \int_0^t Z^i_sZ^j_s\, \ud s}, \qquad j \in D,
	\end{equation}
	with starting value $\VEC{z}=(z^i, i\in D)\in \re_+^d$.

	When $\mat{C}\equiv\mat{0}$, the authors in \cite{CPGUB} proved, by using analytical techniques, the existence and uniqueness of a solution to the multi-parameter time-change equation  \eqref{eq: CIMBP TC}. 
	In the next subsection, we show the existence of a process that is a  solution to \eqref{eq: CIMBP TC}. Such process is obtained by a scaling limit representation. 
	Additionally, we prove that 
	it coincides with $\VEC{Y}$, the unique strong solution to \eqref{continuous MLB}, i.e. it is a CIMBP.
	
	\subsection{Scaling limits of DIMBP}\label{limit}
	The main result of this work is the convergence of re-normalized discrete models to the continuous one.
	\smallskip
	
	Let $(\VEC{X^{(n),i}}, 1\leq i\leq d)_{n \in \mathbb{N}}$ be a sequence of random walks such that $X^{(n),i,i}$ is downwards skip-free and $X^{(n),i,j}$ is non-decreasing for $j\neq i$. Assume that for each $i\in D$, $(X^{(n),i, j}, j\in D)_{n \in \mathbb{N}}$ can be scaled to converge to a L\'evy process $\VEC{X^i}=(X^{i,j}, j\in D)$. That is, assume the existence of  constants $(a_n)_{n \in \mathbb{N}}$ and $(b^i_n, i \in D )_{n \in \mathbb{N}}$   such that 
	\begin{equation} \tag{{\bf CL}}\label{eq: CL}
	\left(\frac{a_n}{b^j_n}X^{(n),i,j}_{b^i_nt },\ t\geq 0 , j \in D \right) \underset{n\to\infty}{\longrightarrow} \left(X^{i,j}_{t},\ t\geq 0, j\in D\right),  \qquad i\in D, 
	\end{equation}
	where $a_n \underset{n\to\infty}{\longrightarrow} \infty$ and $b^j_n/a_n\underset{n\to\infty}{\longrightarrow} \infty$.
	The convergence is almost surely in the Skorokhod space. 
	
	\begin{theorem}\label{thm: SL}
		Assume hypothesis \eqref{eq: CL}. 
		Let $(z^i_n,\ i\in D)_{n \in \mathbb{N}}$ and $(\mat{C_n})_{n \in \mathbb{N}}=(c^{i,j}_n,\ i,j \in D)_{n \in \mathbb{N}}$ be such that 
		$z^j_n a_n/b^j_n \underset{n\to\infty}{\longrightarrow} z^j$  and $c^{i,j}_nb^i_n \underset{n\to\infty}{\longrightarrow}c^{i,j}$. 
		Let   $\VEC{Z^{(n)}}=\left\{\left(Z^{(n),1}_t,\ldots, Z^{(n),d}_t\right),\ t \geq 0\right\}$ be a DIMBP that satisfies equation 
		\eqref{eq: DIMBP TC} for $(\VEC{X^{(n),i}},\ i\in D)$ and $\mat{N}=(N^{i,j},\ i,j \in D)$,  with starting value $\VEC{z}_n=(z^i_n,\ i\in D)$ and competition parameters $\mat{C_n}$. 
		Then, 
		\begin{equation*}
		\paren{\frac{a_n}{b^i_n}Z^{(n),i}_{a_nt},\ t\geq 0, i \in D} \underset{n\to\infty}{\longrightarrow}\paren{Y^i_{t},\ t\geq 0, i \in D }
		\end{equation*}
		almost surely in the Skorokhod space, where $\VEC{Y}=\left\{(Y^{1}_t,\ldots, Y^{d}_t),\ t\geq 0 \right\}$ is the solution to \eqref{continuous MLB}    and it is also a solution of \eqref{eq: CIMBP TC} with starting value $\VEC{z}=(z^i,\ i\in D)$.
	\end{theorem}
	
	To prove this result, first we approximate our interacting processes by a sequence of piece-wise multitype continuous-state branching processes (MCBP) on a time-population grid. Then, we use the results of \cite{CPGUB} for every MCBP to prove the scaling limits for the approximations. Our construction 
	allows us to see that a  CIMBP is also a generalized Lamperti-type  transformation of multidimensional L\'evy processes. With these tools, we conclude the desired convergence.
	
	\section{Proof of Theorem \ref{thm: SL} by using a time-population grid approximation}\label{se: Grid}
	
	We start by constructing a sequence of piecewise multitype continuous-state branching processes on a time-population grid which approximate a CIMBP.  Let $\epsilon, \delta > 0$ and $\VEC{z}=(z^i, i\in D)$. For each $j \in D$, we define $k^{0,j}:= \lfloor\tfrac{z^j}{\delta}\rfloor $. 
	By \cite[Theorem 1]{CPGUB}  the time-change equation
	\begin{equation}\label{eq: z1}
	Z^{1,j}_t= z^j + \sum\limits_{i=1}^d X^{i,j}_{\int_0^t Z^{1,i}_s \ud s} + \sum\limits_{i=1}^d c^{ij}k^{0,j}\delta\int_0^t Z^{1,i}_s \ud s,
	\quad t\geq 0, \quad j \in D,
	\end{equation}
	has a unique solution, which is a MCBP. Therefore, by using \cite[Theorem 4.6]{MR3340375}, $\VEC{Z^1}$ is also the unique strong solution to 
	\begin{eqnarray*}
		Z^{1,j}_t&=&z^j+  \sum\limits_{i=1}^d c^{ij}k^{0,j}\delta\int_0^t Z^{1,i}_s\ud s +\sum\limits_{i=1}^d\int_0^t b^{ij}Z^{1,i}_s\ud s + \int_0^t \sqrt{2\sigma^jZ^{1,j}_s}\ud W^{(j)}_s\\
		&& +\sum\limits_{i=1}^d \int_0^t\int_{\re_+^d}\int_0^{\infty}r_j\indi{\{u\leq Z^{1,i}_{s-}\}} \widetilde{\mathcal{N}}^{(i)}(\ud s,\ud \VEC{r},\ud u), 
		\qquad t\geq 0, \quad j \in D. 
	\end{eqnarray*}
	where $(W^{(i)}, i\in D)$ are standard Brownian motions and $(\mathcal{N}^{(i)},i \in D)$ are Poisson random measures on $\re_+^{d+2}$ 
	defined as in equation \eqref{continuous MLB}.
	
	Given $(\VEC{Z}^{\VEC{1}}_t,\ t\geq 0)$, for each $i \in D$, the vector-valued process defined by
	\begin{equation}\label{X1}
	X^{\{1\},i,j}_r := X^{i,j}_{\int_{0}^\varepsilon Z^{1,i}_s \ud s + r} - X^{i,j}_{\int_{0}^\varepsilon Z^{1,i}_s \ud s}, \qquad j \in D,
	\end{equation}
	is a L\'evy process in $\re^d$. If we set $k^{1,j}:=\lfloor \tfrac {Z^{1,j}_{\varepsilon}} {\delta }\rfloor$, again the time-change equation
	\begin{equation}\label{eq: z2} 
	Z^{2,j}_t=  Z^{1,j}_\varepsilon + \sum\limits_{i=1}^d X^{\{1\},i,j}_{\int_{0}^{t}Z^{2,i}_s \ud s} + \sum\limits_{i=1}^d c^{i,j}\left(k^{1,j}\delta\right)\int_0^t Z^{2,i}_s \ud s,
	\quad t\geq 0, \quad j \in D.
	\end{equation}
	has an  unique solution, which is a MCBP starting from $\VEC{Z}^{\VEC{1}}_\varepsilon$.
	
	Inductively, for each $m\geq 1$, given $(\VEC{Z}^{\VEC{m}}_t,\ t\geq 0)$, we define $k^{m,j}:=\lfloor \tfrac {Z^{m,j}_{\varepsilon}} {\delta }\rfloor$. The time-change equation
	\begin{equation}\label{eq: Zm}
	Z^{m+1,j}_t=  Z^{m,j}_\varepsilon + \sum\limits_{i=1}^d X^{\{m\},i,j}_{\int_{0}^{t}Z^{m+1,i}_s \ud s} + \sum\limits_{i=1}^d c^{i,j}\left(k^{m,j}\delta\right)\int_0^t Z^{m+1,i}_s \ud s,
	\quad t\geq 0, \quad j \in D,
	\end{equation}
	has a unique solution, where for each $r\geq 0$
	\begin{equation}\label{eq: x{m}}
	\begin{split}
	X^{\{m\},i,j}_r &= X^{\{m-1\},i,j}_{\int_{0}^\varepsilon Z^{m,i}_s \ud s + r} - X^{\{m-1\},i,j}_{\int_{0}^\varepsilon Z^{m,i}_s \ud s}\\ 
	&=X^{i,j}_{\left(\sum\limits_{k=1}^m \int_{0}^\varepsilon Z^{k,i}_s \ud s\right) + r} - X^{i,j}_{\sum\limits_{k=1}^m \int_{0}^\varepsilon Z^{k,i}_s \ud s }. 
		\end{split}
	\end{equation}
	
	Observe that $\VEC{Z}^{\VEC{m}}$ is a MCBP starting from $\VEC{Z}^{\VEC{m-1}}_\varepsilon$. Again, by using \cite[Theorem 4.6]{MR3340375}, $\VEC{Z}^{\VEC{m}}$ is also the unique strong solution to
	\begin{eqnarray}
	\label{eq: Ymj}
	\begin{split}
	Z^{m,j}_t=& Z^{m-1,j}_\varepsilon +  \sum\limits_{i=1}^d c^{ij}k^{m-1,j}\delta \int_{0}^{t} Z^{m,i}_s\ud s +\sum\limits_{i=1}^d\int_{0}^{t} b^{ij}Z^{m,i}_s\ud s \\ 
	 &+ \int_0^{t} \sqrt{2\sigma^jZ^{m,j}_s}\ud W^{(m-1,j)}_s
	+\sum\limits_{i=1}^d \int_{0}^{t}\int_{\re_+^d}\int_0^{\infty}r_j\indi{\{u\leq Z^{m,i}_{s-}\}} \widetilde{\mathcal{N}}^{(m-1,i)}(\ud s,\ud \VEC{r},\ud u) ,
		\end{split}
	\end{eqnarray}
	where 
	$
	W^{(m-1,j)}_t=W^{(j)}_{t+(m-1)\varepsilon}-W^{(j)}_{(m-1)\varepsilon}$ is a Brownian motion and the measure $\mathcal{N}^{(m-1,i)}$ is given by 
	$\mathcal{N}^{(m-1,i)}(\ud s,\ud \VEC{r},\ud u )= \mathcal{N}^{(i)}(\ud s + (m-1)\varepsilon,\ud \VEC{r},\ud u ).
	$
	
	Using the processes constructed above,  we define the process $(\VEC{Z}^{\varepsilon,\delta}_t,t\geq 0)$ as
	\begin{equation*}\label{eq: partz}
	\VEC{Z}^{\varepsilon,\delta}_t:= 
	\VEC{Z}^{m+1}_{t - m\varepsilon},\qquad t \in [m\varepsilon,(m+1)\varepsilon], \quad m\in\{0,1,2,\dots\}.
	\end{equation*}
	
	
	The previous process is well defined at multiples of $\varepsilon$ since $\VEC{Z}_{(m+1)\varepsilon}^{\varepsilon,\delta}= \VEC{Z}^{m+1}_{0}=\VEC{Z}^{m}_{\varepsilon}$ for $m\in \mathbb{N}$, where the last identity follows by \eqref{eq: Zm}.
	On the other hand, if $t \in (m\varepsilon,(m+1)\varepsilon)$ using equations \eqref{eq: Zm} and \eqref{eq: x{m}}, together with the definition of $\VEC{Z}^{\varepsilon,\delta}$, we have that
	\begin{eqnarray*}
		Z^{\varepsilon,\delta,j}_t - Z^{\varepsilon,\delta,j}_{m\varepsilon}&=& \sum\limits_{i=1}^d \left[X^{i,j}_{\left(\sum\limits_{k=1}^{m}\int_0^\varepsilon Z^{k,i}_s \ud s\right) +\int_0^{t-m\varepsilon} Z^{m+1,i}_s \ud s } -
		X^{i,j}_{\sum\limits_{k=1}^{m}\int_0^\varepsilon Z^{k,i}_s \ud s}\right]\\ \nonumber
		&&+ \sum\limits_{i=1}^d c^{i,j}\left(\lfloor \tfrac {Z^{m,j}_{\varepsilon}} {\delta }\rfloor\delta\right)\int_0^{t-m\varepsilon} Z^{m+1,i}_s \ud s\\ \nonumber
		&=&\sum\limits_{i=1}^d \left[X^{i,j}_{\int_0^t Z^{\varepsilon,\delta,i}_s \ud s} -
		X^{i,j}_{\int_0^{m\varepsilon} Z^{\varepsilon,\delta,i}_s \ud s}\right]+ \sum\limits_{i=1}^d c^{i,j}\left(\lfloor \tfrac {Z^{\varepsilon,\delta,j}_{m\varepsilon}} {\delta }\rfloor\delta\right)\int_{m\varepsilon}^t Z^{\varepsilon,\delta,i}_s \ud s.
	\end{eqnarray*}
	Thus, the process $\VEC{Z^{\varepsilon,\delta}}$ satisfies 
	\begin{equation}\label{eq: Zepsdelta}
	Z^{\varepsilon,\delta,j}_t=z^j + \sum\limits_{i=1}^d X^{i,j}_{\int_{0}^t Z^{\varepsilon,\delta,i}_s \ud s} +
	\sum\limits_{i=1}^d c^{i,j}\sum\limits_{m=1}^{\lfloor t/\varepsilon \rfloor + 1} \int_{(m-1)\varepsilon}^{m\varepsilon \wedge t} \left\lfloor \frac{Z^{\varepsilon,\delta,j}_{(m-1)\varepsilon}}{ \delta} \right\rfloor\delta Z_s^{\varepsilon,\delta,i}\ud s , \qquad  t\geq 0, j \in D.
	\end{equation}

	If now, we work in a similar way with equation \ref{eq: Ymj}, we deduce that  $Z^{\varepsilon,\delta,j}$ is also the unique solution to
	\begin{eqnarray}\label{eq: Yepsdelta}
	Z^{\varepsilon,\delta,j}_t&=&z^j  + \sum\limits_{i=1}^d c^{i,j}\sum\limits_{m=1}^{\lfloor t/\varepsilon \rfloor + 1}\int_{(m-1)\varepsilon}^{m\varepsilon \wedge t} \left\lfloor \frac{Z^{\varepsilon,\delta,j}_{(m-1)\varepsilon}}{ \delta} \right\rfloor\delta Z_s^{\varepsilon,\delta,i}\ud s +\sum\limits_{i=1}^d\int_{0}^{t} b^{ij}Z^{\varepsilon,\delta,i}_s\ud s \\ \nonumber 
	&& + \int_0^{t} \sqrt{2\sigma^j Z^{\varepsilon,\delta ,j}_s}\ud W^{(j)}_s
	+\sum\limits_{i=1}^d \int_{0}^{t}\int_{\re_+^d}\int_0^{\infty}r_j\indi{\{u\leq Z^{\varepsilon,\delta,i}_{s-}\}} \widetilde{\mathcal{N}}^{(i)}(\ud s,\ud \VEC{r},\ud u).
	\end{eqnarray}
	
	In order to prove Theorem \ref{thm: SL}, we are going to use some auxiliary results, whose proofs can be found in Section \ref{se: proofs}.
	The following proposition states that the above piecewise MCBPs converge to a CIMBP. It also establishes that a  CIMBP can be seen as a generalized Lamperti-type  transformation of a given multidimensional L\'evy process.
	
	\begin{proposition}\label{prop: conv epsdelta}
		Let $(\VEC{Z}^{\VEC{\varepsilon,\delta}}_t, t\geq 0)$ be the unique solution to \eqref{eq: Zepsdelta}, which coincides with the unique solution to \eqref{eq: Yepsdelta}. Then, as $\varepsilon$, $\delta$ go to zero, $\VEC{Z}^{\varepsilon,\delta}$ converges almost surely in the Skorokhod space to $\VEC{Y}$, the unique solution of \eqref{continuous MLB}. Furthermore,  $\VEC{Y}$ satisfies equation \eqref{eq: CIMBP TC}.
	\end{proposition}
	
	Now, for any fixed pair $\varepsilon,\delta> 0$, we prove the scaling limits for the previous approximations. The main idea, as one can see in Section \ref{se: proofs}, is to use the results of \cite{CPGUB} for every piecewise MCBP.
	\begin{lemma}\label{le: SL epsdelta} \textbf{\textsc{Scaling limits on the $\varepsilon,\delta$-grid.}}
		Given $\varepsilon,\delta> 0$, let   $(\VEC{Z}_t^{\varepsilon,\delta,(n)},\ t \geq 0)$ be the process defined as the unique solution to 
		\begin{eqnarray} \label{eq: DIMBP TC epsdelta}
		Z^{\varepsilon,\delta,(n),j}_t&=&z_n^j +\sum_{i=1}^d X^{(n),i,j}_{\int_0^tZ^{\varepsilon,\delta,(n),i}_r\ud r}\\ \nonumber
		&&+\sum_{i=1}^d sgn(c_n^{i,j})N^{i,j}_{|c_n^{i,j}|\sum\limits_{_m=1}^{\lfloor t/a_n\varepsilon \rfloor+1}\int_{ma_n\varepsilon}^{t \wedge (m+1)a_n\varepsilon}\lfloor \frac{Z^{\varepsilon,\delta,(n),j}_{ma_n\varepsilon}}{\delta}\rfloor \delta  Z^{\varepsilon,\delta,(n),i}_r  \ud r},
		\end{eqnarray}
		$j \in D$, driven by $(\VEC{X^{(n),i}}, i\in D)$ and $\mat{N}=(N^{i,j}, i,j \in D)$,  with starting value $\VEC{z}_n=(z_n^j , j\in D)$ and competition parameters $\mat{C_n}$ as in Theorem \ref{thm: SL}.
		Assume hypothesis \eqref{eq: CL}. 
		Then, 
		\begin{equation*}
		\paren{\frac{a_n}{b^i_n}Z^{\varepsilon,\delta,(n),i}_{a_nt},\  t\geq 0, i \in D}\underset{n\to\infty}{\longrightarrow}
		\paren{Z^{\varepsilon,\delta,i}_{t},\ t\geq 0, i \in D }
		\end{equation*}
		almost surely in the Skorokhod space, where $\VEC{Z}^{\varepsilon,\delta}$ is the unique solution to \eqref{eq: Yepsdelta}.
	\end{lemma}
	
	Finally, for each pair $\varepsilon,\delta > 0$ we can write
	\begin{equation*}
	\frac{a_n}{b^j_n}Z^{(n),j}_{a_nt} - Y^{j}_t = \frac{a_n}{b^j_n}\left(Z^{(n),j}_{a_nt} - Z^{\varepsilon,\delta,(n),j}_{a_nt}\right) + \left( \frac{a_n}{b^j_n}Z^{\varepsilon,\delta,(n),j}_{a_nt} - Z^{\varepsilon,\delta,j}_t\right) + \left(Z^{\varepsilon,\delta,j}_t - Y^{j}_t \right),
	\end{equation*}
	for $\VEC{Z}^{\varepsilon,\delta, (n)}$ as in \eqref{eq: DIMBP TC epsdelta} and  $\VEC{Z}^{\varepsilon, \delta}$ the unique solution to \eqref{eq: Zepsdelta}. 
	Thanks to  Proposition \ref{prop: conv epsdelta} and Lemma \ref{le: SL epsdelta}, it only remains to prove that the first term in the right-hand side of the equation converges to zero. This is asserted in the following lemma.
	
	\begin{lemma}\label{le: ZnZn,eps}
		Let $(\VEC{Z}^{\VEC{(n)}}_t,\ t\geq 0)$ be the process referred in Theorem \ref{thm: SL}, and let   $(\VEC{Z}_t^{\varepsilon,\delta,(n)},\ t \geq 0)$ be the process referred in Lemma \ref{le: SL epsdelta}. Then,
		\begin{equation*}
		\lim\limits_{\varepsilon,\delta \rightarrow 0}\lim\limits_{n\rightarrow \infty}\frac{a_n}{b^j_n}\left(Z^{(n),j}_{a_nt} - Z^{\varepsilon,\delta,(n),j}_{a_nt}\right)=0
		\end{equation*}
		almost surely in the  Skorokhod space.
	\end{lemma}
	
	Therefore, from Proposition \ref{prop: conv epsdelta}, Lemma \ref{le: SL epsdelta} and Lemma \ref{le: ZnZn,eps}, we conclude that the re-normalized DIMBPs converge to a CIMBP, i.e. 
	\begin{equation*}
	\paren{\frac{a_n}{b^i_n}Z^{i,(n)}_{a_nt},\ t\geq 0, i \in D}\underset{n\to\infty}{\longrightarrow}
	\paren{Y^i_{t},\ t\geq 0, i \in D }.
	\end{equation*}
	\hfill \qedsymbol
	
	\section{Other proofs}\label{se: proofs}
	In this section, we prove Theorem \ref{thm: continuous MLB},  Proposition \ref{prop: conv epsdelta}, Lemma \ref{le: SL epsdelta} and Lemma \ref{le: ZnZn,eps}.

	\subsection*{Proof of Theorem \ref{thm: continuous MLB}.}
	Following the same ideas used in the proof of \cite[Proposition 2.4]{MR2584896}, the existence of a unique strong solution to
	\begin{equation}\label{eq: Ym}
	\begin{split}
	Y_t^j&=y^j+ \sum\limits_{i=1}^d\int_0^t c^{ij}(Y^i_s \wedge n)(Y^j_s \wedge n)\ud s + \sum\limits_{i=1}^d\int_0^t b^{ij}(Y^i_s \wedge n) \ud s\\ 
	&+\int_0^t \sqrt{2\sigma^j\left(Y^j_s \wedge n\right)}\ud W^{(j)}_s +\sum\limits_{i=1}^d \int_0^t\int_{\re_+^d}\int_0^{\infty}(r_j\wedge n)\indi{\{u\leq (Y^i_{s-}\wedge n)\}} \widetilde{\mathcal{N}}^{(i)}(\ud s,\ud \VEC{r},\ud u), \quad j \in D,
	\end{split}
	\end{equation}
	for each $n\geq 1$, implies the existence of a unique strong solution to \eqref{continuous MLB}.
	
	By \cite[Example 134]{SITU2005}, we know that there exists a unique weak solution for equation \eqref{eq: Ym}. Therefore, by \cite[Theorem 1.2]{MR3298537}, it is enough to show pathwise uniqueness to assure the existence of a unique strong solution. 
	To show the pathwise uniqueness for equation \eqref{eq: Ym}, we will adapt the proof of \cite[Theorem 2.1]{MR2952093}. 
	
	Let $a_0 = 1$ and choose a decreasing sequence  $a_k \rightarrow  0$ such that 
	$\int_{a_k}^{a_{k-1}} z^{-1} \ud z = k$ for $k \geq 1$. Let $x \mapsto \psi_k (x)$ be a non-negative continuous function on $\re$ with support in $(a_k , a_{k-1})$ that satisfies $\int_{a_k}^{a_{k-1}} \psi_k(x)\ud x=1$ and $0\leq \psi_k(x)\leq 2k^{-1}x^{-1}$ for  $a_k<x<a_{k-1}$.  
	For each $k \geq 1$, we define the non-negative and twice continuously differentiable function
	\begin{equation*}
	\varphi_k(z)=\int_0^{|z|}\ud y\int_0^y \psi_k(x)\ud x.
	\end{equation*}
	We observe that the sequence of functions $\{\varphi_k, k\geq 1\}$ satisfies the following properties:
	\begin{itemize}
		\item[(i)]$\varphi_k (z) \rightarrow |z|$ non-decreasingly as $k \rightarrow \infty$;
		\item[(ii)] $ 0 \leq \varphi'_k (z) \leq 1$  for $z \geq 0$ and $-1 \leq \varphi'_k (z) \leq 0$ for $z \leq 0$;
		\item[(iii)] $\varphi''_k(z)\geq 0$ for every $ z \in \re$. In addition, as $k \rightarrow \infty$,
		\[
		\varphi''_k(x-y)(\sqrt{x}-\sqrt{y})^2 \rightarrow 0.
		\]
	\end{itemize}
	
	We first fix the integer $n \geq 1$.  Let $\VEC{Y}=\left\{(Y^1_t,\dots,Y^d_t),\ t\geq 0\right\}$ and $\widetilde{\VEC{Y}}=\left\{(Y^1_t,\dots,Y^d_t),\ t\geq 0 \right\}$  be two  weak solutions of \eqref{eq: Ym}.	For each $j \in D$, we define $\zeta^j_t=Y^j_t - \widetilde{Y}^j_t$ for $t \geq 0$. Then, by It\^o’s formula, for $j\in D$
	\begin{equation}
	\label{ineq1}
	\begin{split}
	\varphi_k(\zeta^j_{t})= &\sum\limits_{i=1}^d\int_0^{t} \varphi'_k(\zeta^j_{s-}) \left(b^{ij} + c^{ij}(\widetilde{Y}^j \wedge n)\right)\left(Y^i_s\wedge n - \widetilde{Y}^i_s\wedge n\right)\ud s\\ 
	&+ \sum\limits_{i=1}^d\int_0^{t} \varphi'_k(\zeta^j_{s-})c^{ij}(Y^i_s \wedge n)\left(Y^j_s \wedge n - \widetilde{Y}^j_s \wedge n\right) \ud s\\ 
	& + \int_0^{t} \varphi'_k(\zeta^j_{s-})\sqrt{2\sigma^j}\left[\sqrt{Y^j_s\wedge n} -\sqrt{\widetilde{Y}^j_s\wedge n}\right]\ud W^{(j)}_s \\
	& + \int_0^{t} \varphi''_k(\zeta^j_{s-})\sigma^j\left[\sqrt{Y^j_s\wedge n} -\sqrt{\widetilde{Y}^j_s\wedge n}\right]^2 \ud s \\ 
	& + \sum\limits_{i=1}^d\int_0^{t}\int_{\re_+^d}\int_0^{\infty}\left[\indi{\{u\leq Y^i_{s-}\wedge n\}}-\indi{\{u\leq \widetilde{Y}^i_{s-}\wedge n\}} \right]\Delta_{r_j \wedge n}\varphi_k(\zeta^j_{s-}) \widetilde{\mathcal{N}}^{(i)}(\ud s,\ud \VEC{r},\ud u)\\ 
	&+\sum\limits_{i=1}^d \int_0^{t}\int_{\re_+^d}\left(Y^i_{s-}\wedge n - \widetilde{Y}^i_{s-}\wedge n\right)D_{r_j \wedge n}\varphi_k(\zeta^j_{s-}) \ud s m^{(i)}(\ud \VEC{r})
		\end{split}
	\end{equation}
	where $\Delta_{r}\varphi(z)=\varphi(z + r)- \varphi(z)$ and $D_r\varphi(z)= \Delta_r\varphi(z) - \varphi'(z)r$. We want to ensure that the integral with respect to  $m^{(i)}$ is bounded. We recall that the measure $m^{(i)}$  satisfies the integrability condition \eqref{eq: integral conditions}. Therefore, we need to treat $\{j=i\}$ separately. Note that 
	\begin{equation*}
\begin{split}
		\left(Y^i_{s-}\wedge n - \widetilde{Y}^i_{s-}\wedge n\right)&D_{r_j \wedge n}\varphi_k(\zeta^j_{s-})\\
	=	&\left(Y^i_{s-}\wedge n - \widetilde{Y}^i_{s-}\wedge n\right)(r_j \wedge n)^2 \int_0^1 \psi_k(|\zeta^j_{s-}+t(r_j \wedge n)|)(1-t)\ud t\\ \nonumber
	\leq	&  \left(Y^i_{s-}\wedge n - \widetilde{Y}^i_{s-}\wedge n\right)(r_j \wedge n)^2 \int_0^1 \dfrac{2(1-t)}{k|\zeta^j_{s-}+t(r_j \wedge n)|}\ud t \\ \nonumber
	\leq	& \left\{
		\begin{array}{rr}
			\frac{(r_j \wedge n)^2}{k} & \quad \text{ if }\quad \zeta^j_{s-}>0\\
			0       & \quad \text{ if }\quad \zeta^j_{s-}\leq 0
		\end{array}
		\right.,
		\end{split}
	\end{equation*}
	and also  
	\begin{equation*}
	D_{r_j \wedge n}\varphi_k(\zeta^j_{s-}) = \Delta_{r_j \wedge n}\varphi_k(\zeta^j_{s-}) - \varphi_k'(\zeta^j_{s-})(r_j \wedge n) \leq 2(r_j \wedge n).
	\end{equation*}
	
	We apply the previous bounds to equation \eqref{ineq1} to get
	\begin{eqnarray*}
		\varphi_k(\zeta^j_{t})	
		&\leq & \sum\limits_{i=1}^d \int_0^{t} \varphi'_k(\zeta^j_{s-})\left(b^{ij} + c^{ij}(\widetilde{Y}^j \wedge n)\right)\left(Y^i_s\wedge n - \widetilde{Y}^i_s\wedge n\right)\ud s\\ \nonumber
		&& +  \sum\limits_{i=1}^d \int_0^{t} \varphi'_k(\zeta^j_{s-})c^{ij}(Y^i_s \wedge n)\left(Y^j_s \wedge n -\widetilde{Y}^j_s \wedge n\right) \ud s\\\nonumber
		&& +\int_0^{t} \varphi''_k(\zeta^j_{s-})\sigma^j\left[\sqrt{Y^j_s\wedge n} -\sqrt{\widetilde{Y}^j_s\wedge n}\right]^2 \ud s + \textit{local martingale}	\\ \nonumber
		&& + \sum\limits_{i=1,i\neq j}^d \int_0^{t}\left(Y^i_{s-}\wedge n - \widetilde{Y}^i_{s-}\wedge n\right)\int_{\re_+^d}2(r_j \wedge n) m^{(i)}(\ud \VEC{r}) \ud s\\ \nonumber
		&& + \int_0^{t}\left(Y^j_{s-}\wedge n - \widetilde{Y}^j_{s-}\wedge n\right)\int_{\re_+^d}2(r_j \wedge n)\indi{\{r_j\geq 1\}} m^{(j)}(\ud \VEC{r})\ud s\\ \nonumber
		&& + \int_0^{t}\int_{\re_+^d}\dfrac{(r_j \wedge n)^2}{k}\indi{\{r_j\leq 1\}}  m^{(j)}(\ud \VEC{r})\ud s.
	\end{eqnarray*}
	If we take expectation and let $k \rightarrow \infty$, by the Monotone Convergence Theorem, we  obtain
	\begin{eqnarray*}
		\se\left(|\zeta^j_{t}|\right)	
		&\leq & \sum\limits_{i=1}^d \int_0^t \left[|b^{ij}+ c^{ij}n|\se\left(|Y^i_s\wedge n - \widetilde{Y}^i_s\wedge n|\right) + |c^{ij}n|\se\left(|Y^j_s \wedge n -\widetilde{Y}^j_s \wedge n|\right) \right] \ud s\\ \nonumber
		&& + 2n \sum\limits_{i=1,i\neq j}^d \left[\int_{\re_+^d} m^{(i)}(\ud \VEC{r})\right]\int_0^t\se\left(|Y^i_s\wedge n - \widetilde{Y}^i_s\wedge n|\right)\ud s \\ \nonumber
		&&+ 2n\left[\int_{\re_+^d}\indi{\{r_j\geq 1\}} m^{(j)}(\ud \VEC{r})\right]\int_0^t\se\left(|Y^j_s \wedge n -\widetilde{Y}^j_s \wedge n|\right)\ud s\\ \nonumber
		&\leq& C_j\sum\limits_{i=1}^d \int_0^t \se\left(|Y^i_s\wedge n - \widetilde{Y}^i_s\wedge n|\right)\ud s
		\leq C_j\sum\limits_{i=1}^d \int_0^t \se\left(|\zeta^i_{s}| \wedge n \right)\ud s\\ \nonumber
		&\leq &C_j\sum\limits_{i=1}^d \int_0^t \se\left(|\zeta^i_{s}| \right)\ud s
	\end{eqnarray*}
	for a big enough constant $C_j$. Thus, by summing over $j$, we have 
	\begin{equation*}
	\se\left(\sum\limits_{j=1}^d|\zeta^j_{t}|\right) \leq \bar{C} \int_0^t \se\left(\sum\limits_{j=1}^d|\zeta^j_{s}|\right)\ud s,
	\end{equation*}
	where $\bar{C}=\sum\limits_{j=1}^dC_j$. Therefore,  by Gronwall's inequality,
	\begin{equation*}
	\se\left(\sum\limits_{j=1}^d|\zeta^j_{t}|\right)=0 \qquad \text{ for all } t\geq 0.
	\end{equation*}
	This implies the pathwise uniqueness for equation \eqref{eq: Ym} for all $n \geq 1$, and from here we deduce that there exists a unique strong solution for \eqref{continuous MLB}. 
	
	\smallskip
	
	Finally, by It\^o's formula for every $\VEC{F} \in C^2 (\re^d_+)$ we have 
	\begin{eqnarray*}
		\VEC{F}(\VEC{Y}_t) &=& \VEC{F}(\VEC{Y}_0)+\int_0^t \sum\limits_{i=1}^d\left\{F''_{ii}(\VEC{Y}_{s-})\sigma_iY^i_s + F'_i(\VEC{Y}_{s-})\sum\limits_{j=1}^d\left[b^{ij}Y^j_s + c^{ij}Y^j_sY^j_s\right] \right\}\ud s\\ \nonumber
		&& +  \sum\limits_{i=1}^d \int_0^t\int_{\re_+^d} Y^i_s \left[\VEC{F}(\VEC{Y}_{s-}+\VEC{r}) - \VEC{F}(\VEC{Y}_{s-}) -\sum\limits_{j=1}^dF'_j(\VEC{Y}_{s-})r_j\right]m^{(i)}(\ud \VEC{r})\ud s + M_t,
	\end{eqnarray*}
	where $M_t$ is a local martingale. Then, $\VEC{Y}$ has the desired infinitesimal generator.
	\hfill \qedsymbol
	
	\subsection*{Proof of Proposition \ref{prop: conv epsdelta}.}
	As some of the techniques are the same as those used in the previous proof, we only give the short version.
	Let $\VEC{Y}= \left\{(Y^1_t,\dots,Y^d_t),\ t\geq 0\right\}$ be the unique strong solution of \eqref{continuous MLB} and for every $\varepsilon,\delta>0$ let  $\VEC{Z}^{\varepsilon,\delta}=\left\{(Z^{\varepsilon,\delta,1}_t,\dots,Z^{\varepsilon,\delta,d}_t),\ t \geq 0\right\}$  be the strong solution of \eqref{eq: Yepsdelta}. For each $j \in \{1,\dots,d\}$, we define $\Delta^j_t=\Delta^{\varepsilon,\delta,j}_t=Y^j_t - Z^{\varepsilon,\delta,j}_t$.
	
	Let $\tau_M:=\inf\{t\geq 0: Y^i_t>M \mbox{ or } Z^{\varepsilon,\delta,i}_t>M \mbox{ for some } i \in D\}$. Take a sequence of functions $\{\varphi_k\}_{k\in \mathbb{N} }$ as in the proof of Theorem \ref{thm: continuous MLB}. By 
	Itô’s formula, we have 
	\begin{eqnarray*}
		\varphi_k(\Delta^j_{t \wedge \tau_M})&=&\sum\limits_{i=1}^d\int_0^{t\wedge \tau_M} \varphi'_k(\Delta^j_{s-}) \left(b^{ij} + c^{ij}Y_s^j\right)\Delta^i_s\ud s
		+ \sum\limits_{i=1}^d\int_0^{t \wedge \tau_M} \varphi'_k(\Delta^j_{s-}) c^{ij}Z^{\varepsilon,\delta,i}_s\Delta^j_s\ud s\\ \nonumber
		&&+\sum\limits_{i=1}^dc^{i,j}\sum\limits_{m=1}^{\lfloor (t\wedge \tau_M)/\varepsilon \rfloor + 1}\int_{(m-1)\varepsilon}^{m\varepsilon \wedge (t\wedge \tau_M)} \varphi'_k(\Delta^j_{s-})\left[Z^{\varepsilon,\delta,j}_s - \left\lfloor \frac{Z^{\varepsilon,\delta,j}_{(m-1)\varepsilon}}{ \delta} \right\rfloor\delta \right]Z_s^{\varepsilon,\delta,i}\ud s\\ \nonumber
		&& + \int_0^{t\wedge \tau_M} \varphi''_k(\Delta^j_{s-})\sigma^j\left[\sqrt{Y^j_s} -\sqrt{Z^{\varepsilon,\delta,j}_s}\right]^2 \ud s + \textit{local martingale}	\\ \nonumber
		&&+\sum\limits_{i=1}^d \int_0^{t\wedge \tau_M}\int_{\re_+^d}\Delta^i_sD_{r_j}\varphi_k(\Delta^j_{s-}) \ud s m^{(i)}(\ud \VEC{r})
	\end{eqnarray*}
	
	By using the properties of $\{\varphi_k\}_{k\in \mathbb{N} }$,  taking expectation and letting $k \rightarrow \infty$, we obtain 
	\begin{equation*}
		\begin{split}
		\se\left(|\Delta^j_{t\wedge \tau_M}|\right)	
		&\leq \sum\limits_{i=1}^d \int_0^{t\wedge \tau_M} \left[|b^{ij}+ c^{ij}M|\se\left(|\Delta^i_s|\right) + |c^{ij}|M\se\left(|\Delta^j_s|\right) \right] \ud s\\ \nonumber
		+&dM\sum\limits_{m=1}^{\lfloor t/\varepsilon \rfloor + 1}\int_{(m-1)\varepsilon}^{m\varepsilon \wedge t}\se\left(|Z^{\varepsilon,\delta,j}_{s\wedge \tau_M} - \left\lfloor \frac{Z^{\varepsilon,\delta,j}_{(m-1)\varepsilon\wedge \tau_M}}{ \delta} \right\rfloor\delta |\right)\ud s\\ \nonumber
		 +& 2M\left[ \int_{\re_+^d}\indi{\{r_j\geq 1\}} m^{(j)}(\ud \VEC{r})\int_0^{t\wedge \tau_M}\se\left(|\Delta^j_s|\right)\ud s+\sum\limits_{i=1,i\neq j}^d \int_{\re_+^d} m^{(i)}(\ud \VEC{r})\int_0^{t\wedge \tau_M}\se\left(|\Delta^i_s|\right)\ud s\right] \\ \nonumber
		\leq& C_j\sum\limits_{i=1}^d \int_0^{t} \se\left(|\Delta^i_{s\wedge \tau_M}|\right)\ud s 
		+ dM \sum\limits_{m=1}^{\lfloor t/\varepsilon \rfloor + 1}\int_{(m-1)\varepsilon}^{m\varepsilon \wedge t}\se\left(|Z^{\varepsilon,\delta,j}_{s\wedge \tau_M} - \left\lfloor \frac{Z^{\varepsilon,\delta,j}_{(m-1)\varepsilon\wedge \tau_M}}{ \delta} \right\rfloor\delta |\right)\ud s
			\end{split}
	\end{equation*}
	for a big enough constant $C_j$. 
	Thus we have 
	\begin{equation*}
	\begin{split}
		\se\left(\sum\limits_{j=1}^d|\Delta^j_{t\wedge \tau_M}|\right) &\leq  \bar{C} \int_0^{t} \se\left(\sum\limits_{j=1}^d|\Delta^j_{s\wedge \tau_M}|\right)\ud s
		+ dM\sum\limits_{j=1}^d\sum\limits_{m=1}^{\lfloor t/\varepsilon \rfloor + 1}\int_{(m-1)\varepsilon}^{m\varepsilon \wedge t} \se \left(|Z^{\varepsilon,\delta,j}_{s\wedge \tau_M} -Z^{\varepsilon,\delta,j}_{(m-1)\varepsilon\wedge\tau_M} |\right)\ud s\\ \nonumber
		&+ dM\sum\limits_{j=1}^d\sum\limits_{m=1}^{\lfloor t/\varepsilon \rfloor + 1} \varepsilon\,\se\left(\left|Z^{\varepsilon,\delta,j}_{(m-1)\varepsilon\wedge\tau_M} - \left\lfloor \frac{Z^{\varepsilon,\delta,j}_{(m-1)\varepsilon\wedge\tau_M}}{ \delta} \right\rfloor\delta \right|\right),
		\end{split}
	\end{equation*}
	where $\bar{C}=\sum\limits_{j=1}^dC_j$. 
	By Gronwall's inequality, we have that
	\begin{eqnarray*}
		\se\left(\sum\limits_{j=1}^d|\Delta^j_{t\wedge \tau_M}|\right)
		&\leq& e^{\bar{C}t }
		dM\sum\limits_{j=1}^d\sum\limits_{m=1}^{\lfloor t/\varepsilon \rfloor + 1}\int_{(m-1)\varepsilon}^{m\varepsilon \wedge t} \se \left(\left|Z^{\varepsilon,\delta,j}_{s \wedge\tau_M} -Z^{\varepsilon,\delta,j}_{(m-1)\varepsilon\wedge\tau_M} \right|\right)\ud s \\ \nonumber 
		&& +   \varepsilon \sum\limits_{j=1}^d\sum\limits_{m=1}^{\lfloor t/\varepsilon \rfloor + 1} \se\left(\left|Z^{\varepsilon,\delta,j}_{(m-1)\varepsilon\wedge\tau_M} - \left\lfloor \frac{Z^{\varepsilon,\delta,j}_{(m-1)\varepsilon\wedge\tau_M}}{ \delta} \right\rfloor\delta \right|\right).
	\end{eqnarray*}
	We observe that $ |x-\delta \lfloor x/\delta\rfloor|$ goes to zero as $\delta\rightarrow 0$, uniformly in $\re$. 
	
	Additionally, $\left|Z^{\varepsilon,\delta,j}_{s \wedge\tau_M} -Z^{\varepsilon,\delta,j}_{(m-1)\varepsilon\wedge\tau_M} \right|\leq 2M$, and $Z^{\varepsilon,\delta,j}$ is a c\`adl\`ag process. Therefore, by applying the Dominated Convergence Theorem to the right-hand side terms,  we deduce
	\begin{equation*}
	\lim\limits_{\varepsilon,\delta \rightarrow 0}\se\left(\sum\limits_{j=1}^d|\Delta^j_{t \wedge \tau_M}|\right)=0, \qquad \text{ for all } t\geq 0 \mbox{ and } M\geq 0.
	\end{equation*}
	By taking the limit as $M\rightarrow \infty$, we obtain the convergence in the sense of finite dimensional distributions. In order to obtain Skorohod convergence in distribution, by \cite[Theorem 23.9]{KA2021}, it is enough to prove tightness for any sequence  $\{\VEC{Z}^{k}:=\VEC{Z}^{\varepsilon_k, \delta_k}\}_{ k\geq 0}$ when $(\varepsilon_k, \delta_k)\rightarrow 0$. According to Aldous' criterion \cite[Theorem 23.11]{KA2021}, it is enough to prove that for any bounded $\VEC{Z}^k$ - optional times $T_k$ and positive constants $h_k\rightarrow 0$, 
	\begin{equation}\label{eq: convp}
	|Z^{k,j}_{T_k + h_k} - Z^{k,j}_{T_k} |\longrightarrow 0\qquad  \mbox{ as } k\rightarrow \infty \quad \mbox{ in probability}.
	\end{equation}
	Without loosing generality, we can suppose that 
	$h_k\leq  \varepsilon_k/2$. Let  $\tau_M:=\inf\{t\geq 0:  Z^{k,i}_t>M \mbox{ for some } i \in D\}$ and take $T_{k,M}=T_k\wedge \tau_M$ and $m_k=\lfloor \tfrac{T_{k,M}}{\varepsilon_k}\rfloor$. 
	
	Then,  by construction $\lfloor \tfrac{T_{k,M}+h_k}{\varepsilon_k} \rfloor=m_k$ or $\lfloor \tfrac{T_{k,M}+h_k}{\varepsilon_k} \rfloor=m_k+1$.
	By \eqref{eq: Yepsdelta}, we have 
	\begin{eqnarray*}
		|Z^{k,j}_{T_{k,M} + h_k} - Z^{k,j}_{T_{k,M}} |
		&=& \left| \sum\limits_{i=1}^d c^{i,j} \int_{T_{k,M}}^{(m_k+1)\varepsilon_k \wedge (T_{k,M}+h_k)} \left\lfloor \frac{Z^{k,j}_{m_k\varepsilon_k}}{\delta_k}\right\rfloor\delta_k Z^{k,i}_s ds \right.\\
		&& \left. + \sum\limits_{i=1}^d c^{i,j}\int_{(m_k +1)\varepsilon_k \wedge (T_{k,M}+h_k)}^{T_{k,M}+h_k} \left\lfloor \frac{Z^{k,j}_{(m_k+1)\varepsilon_k}}{\delta_k}\right\rfloor\delta_k Z^{k,i}_s ds    \right. \\
		&& +\sum\limits_{i=1}^d b^{ij}\int_{T_{k,M}}^{T_{k,M}+h_k} Z^{k,i}_s\ud s + \int_{T_{k,M}}^{T_{k,M}+h_k} \sqrt{2\sigma^j Z^{k ,j}_s}\ud W^{(j)}_s\\
		&&
		\left.+\sum\limits_{i=1}^d \int_{T_{k,M}}^{T_{k,M}+h_k}\int_{\re_+^d}\int_0^{\infty}r_j\indi{\{u\leq Z^{k,i}_{s-}\}} \widetilde{\mathcal{N}}^{(i)}(\ud s,\ud \VEC{r},\ud u)\right|\\
		&\leq & M h_k \sum\limits_{i=1}^d |c^{i,j}|  \left\lfloor \frac{M}{\delta_k}\right\rfloor\delta_k  
		+Mh_k\sum\limits_{i=1}^d |b^{ij}| + \left|\int_{T_{k,M}}^{T_{k,M}+h_k} \sqrt{2\sigma^j Z^{k ,j}_s}\ud W^{(j)}_s \right|\\
		&&+\sum\limits_{i=1}^d\left| \int_{T_{k,M}}^{T_{k,M}+h_k}\int_{\re_+^d}\int_0^{\infty}r_j\indi{\{u\leq Z^{k,i}_{s-}\}} \widetilde{\mathcal{N}}^{(i)}(\ud s,\ud \VEC{r},\ud u)\right|
	\end{eqnarray*}
	Note that the last two terms correspond to the absolute values of some martingales so, by using Burkholder-Davis-Gundy inequality \cite[Theorem 4.48]{Protter} we have that
	\begin{eqnarray*}
		\mathbb{E}\left[|Z^{k,j}_{T_{k,M} + h_k} - Z^{k,j}_{T_{k,M}} |\right]
		&\leq & M h_k \sum\limits_{i=1}^d |c^{i,j}|  \left\lfloor \frac{M}{\delta_k}\right\rfloor\delta_k  
		+Mh_k\sum\limits_{i=1}^d |b^{ij}| + \sqrt{2\sigma^jMh_k }\\
		&&+\sum\limits_{i=1}^d\sqrt{Mh_k \int_{\re_+^d}r^2_j m^i(\ud\VEC{r})}
	\end{eqnarray*} 
	which goes to zero as $k$ goes to infinity for every $M \in \mathbb{R}_+$. As the $\VEC{Z}^k$ - optional times $T_k$ are bounded, taking $M \rightarrow \infty$ we obtain
	\begin{equation*}
	\lim\limits_{k \rightarrow \infty}\mathbb{E}\left[|Z^{k,j}_{T_{k} + h_k} - Z^{k,j}_{T_{k}} |\right]=0.
	\end{equation*}
	Finally, by Markov's inequality we deduce the desired convergence in probability \eqref{eq: convp}, and from here $\VEC{Z}^{\varepsilon,\delta}$ 
	convergence in distribution in the sense of Skorohod. 
	By the strong uniqueness of the processes (as solutions of SDEs), we have that $ (\VEC{Z}^{\varepsilon,\delta}_t, t\geq 0)$ converges almost surely to $(\VEC{Y}_t, t\geq 0)$ the unique solution to \eqref{continuous MLB}. 
	This implies that for each $t\geq 0$
	\begin{equation*}
	\int_{0}^t Z^{\varepsilon,\delta,i}_s \ud s \overset{a.s.}{\longrightarrow}\int_{0}^t Y^{i}_s \ud s \qquad \text{and}\qquad \sum\limits_{m=1}^{\lfloor t/\varepsilon \rfloor + 1} \int_{(m-1)\varepsilon}^{m\varepsilon \wedge t} \left\lfloor \frac{Z^{\varepsilon,\delta,j}_{(m-1)\varepsilon}}{ \delta} \right\rfloor\delta Z_s^{\varepsilon,\delta,i}\ud s\overset{a.s.}{\longrightarrow}\int_0^t Y^j_sY^i_s \ud s.
	\end{equation*}
	Additionally,  $(\VEC{Z}^{\varepsilon,\delta}_t, t\geq 0)$ is a solution to the time-changed equation \eqref{eq: Zepsdelta}.
	By putting  these two pieces together, we get  that  $(\VEC{Y}_t, t\geq 0)$ is a solution to \eqref{eq: CIMBP TC} for the L\'evy process $\VEC{X}$.
	\hfill \qedsymbol	
	
	\subsection*{Proof of Lemma \ref{le: SL epsdelta}.}  
	First, we observe  that by using a recursive procedure, 
	$\VEC{Z}^{\varepsilon,\delta,(n)}$ can be constructed in an analogous way to $\VEC{Z}^{\varepsilon,\delta}$. 
	Let $k^{0,(n),j}:=\lfloor z^j_n/\delta\rfloor$. We define $\VEC{Z}^{1,(n)}$ as the unique solution of
	\begin{eqnarray*}
		Z^{1,(n),j}_t&=& z^j_n + \sum\limits_{i=1}^d X^{(n),i,j}_{\int_0^t Z^{1,(n),i}_s \ud s} 
		+\sum_{i=1}^d sgn(c^{i,j}_n)N^{i,j}_{|c^{i,j}_n|k^{0,(n),j}\delta\int_0^t Z^{1,(n),i}_r  \ud r}\\ \nonumber 
		&=& z^j_n + \sum\limits_{i=1}^d \tilde{X}^{\{0\},(n),i,j}_{\int_0^t Z^{1,(n),i}_s \ud s},
		\qquad j \in D,
	\end{eqnarray*}
	where $(\VEC{\tilde{X}}^{\{0\},(n),i}, i\in D)$ are independent $d$-dimensional random walks, defined by
	$$\tilde{X}^{\{0\},(n),i,j}_t= X^{(n),i,j}_t + sgn(c^{i,j}_n)N^{i,j}_{|c^{i,j}_n|k^{0,(n),j}\delta t}.$$
	
	Observe that  by hypothesis \eqref {eq: CL},  $\tfrac{a_n}{b^j_n}k^{0,(n),j} \rightarrow k^{0,j}=\lfloor z^j/\delta \rfloor$
	and  $c^{i,j}_nb^i_n \underset{n\to\infty}{\longrightarrow}c^{i,j}$. Therefore,
	\begin{equation*}
	\left(\frac{a_n}{b^j_n}\tilde{X}^{\{0\},(n),i,j}_{b^i_n t},\ t\geq 0, i,j \in D\right) \underset{n\rightarrow \infty}{\longrightarrow}
	\left(X^{i,j}_t + c^{ij}k^{0,j}\delta t,\ t\geq 0, i,j \in D  \right)
	\end{equation*} in the Skorokhod space.
	By Corollary 1 in \cite{CPGUB}, we can see that 
	\begin{equation}\label{eq: convz1}
	\left(\frac{a_n}{b^i_n}Z^{1,(n),i}_{a_nt},\ t\geq 0, i\in D\right) \underset{n\rightarrow\infty}{\longrightarrow}
	\left(Z^{1,i}_t,\ t\geq 0, i \in D \right)
	\end{equation}
	almost surely	in the Skorokhod space, where $\VEC{Z}^{1}$ is the unique solution of \eqref{eq: z1}.
	
	Now, given $\VEC{Z}^{1,(n)}$, the vector-valued process defined by
	\begin{equation*}
	X^{\{1\},(n),i,j}_r := X^{(n),i,j}_{\int_{0}^{a_n\varepsilon} Z^{1,(n),i}_s \ud s + r} - X^{(n),i,j}_{\int_{0}^{a_n\varepsilon} Z^{1,(n),i}_s \ud s}, \qquad j \in D
	\end{equation*}
	is a random walk in $\re^d$ for each $i\in D$, and
	\begin{equation*}
	N^{\{1\},i,j}_r:=N^{i,j}_{|c^{i,j}_n|k^{0,(n),j}\delta\int_{0}^{a_n\varepsilon} Z^{1,(n),i}_s \ud s + r} -N^{i,j}_{|c^{i,j}_n|k^{0,(n),j}\delta\int_{0}^{a_n\varepsilon} Z^{1,(n),i}_s \ud s } , \qquad i,j,\in D
	\end{equation*}
	is a Poisson process.  
	We set $k^{1,(n),j}:=\lfloor \tfrac {Z^{1,(n),j}_{a_n\varepsilon}} {\delta }\rfloor$,  and we define $\VEC{Z}^{2,(n)}$ as the unique solution of 
	\begin{equation*}
	\begin{split}
	Z^{2,(n),j}_t&=  Z^{1,(n),j}_{a_n\varepsilon} + \sum\limits_{i=1}^d X^{\{1\},(n),i,j}_{\int_{0}^{t}Z^{2,(n),i}_s \ud s} 
	+ \sum_{i=1}^d sgn(c^{i,j}_n)N^{\{1\},i,j}_{|c^{i,j}_n|k^{1,(n),j}\delta\int_0^t Z^{2,(n),i}_r  \ud r}\\
	&= Z^{1,(n),j}_{a_n\varepsilon} + \sum\limits_{i=1}^d \tilde{X}^{\{1\},(n),i,j}_{\int_{0}^{t}Z^{2,(n),i}_s \ud s},
	\end{split}
	\quad t\geq 0, \quad j \in D,
	\end{equation*}
	where 
	$$\tilde{X}^{\{1\},(n),i,j}_t= X^{\{1\},(n),i,j}_t + sgn(c^{i,j}_n)N^{\{1\},i,j}_{|c^{i,j}_n|k^{1,(n),j}\delta t}.$$
	By taking $t=\varepsilon$ in  \eqref{eq: convz1}
	\begin{equation*}
	\left(\frac{a_n}{b^i_n}Z^{1,(n),j}_{a_n\varepsilon}, \ j \in D\right)\underset{n\rightarrow \infty}{\longrightarrow}
	\left(Z^{1,j}_{\varepsilon},\  j\in D \right).
	\end{equation*}
	Additionally,
	\begin{equation*}
	\left(\frac{a_n}{b^j_n}\tilde{X}^{\{1\}, (n),i,j}_{b^i_n t},\ t\geq 0, i,j \in D\right) \underset{n\rightarrow \infty}{\longrightarrow}
	\left(X^{\{1\},i,j}_t + c^{ij}k^{1,j}\delta t ,\ t\geq 0, i,j \in D  \right)
	\end{equation*}
	in the Skorokhod space, where $(\VEC{X}_t^{\{1\},i},t \geq 0)$ is the process defined in \eqref{X1} and $k^{1,j}=\lfloor \tfrac {Z^{1,j}_{\varepsilon}} {\delta }\rfloor$. Using again Corollary 1 in \cite{CPGUB}, we have that
	\begin{equation*}
	\left(\frac{a_n}{b^i_n}\VEC{Z}^{2,(n),i}_{a_nt},\ t\geq 0, i\in D\right) \underset{n\rightarrow\infty}{\longrightarrow}
	\left(\VEC{Z}^{2,i}_t,\ t\geq 0, i \in D \right), 
	\end{equation*}
	where $\VEC{Z}^{2,i}$ is the unique solution of \eqref{eq: z2}.

	Inductively, for each $n> 2$, given $(\VEC{Z}^{m-1,(n)}_t, t\geq 0)$ with $k^{m-1,(n),j}=\lfloor \tfrac {Z^{m-1,(n),j}_{a_n\varepsilon}} {\delta }\rfloor$, we define $\VEC{Z}^{m,(n)}$ as the unique solution of
	\begin{eqnarray*}
		Z^{m,(n),j}_t&=&  Z^{m-1,(n),j}_{a_n\varepsilon} + \sum\limits_{i=1}^d X^{\{m-1\},(n),i,j}_{\int_{0}^{t}Z^{m,(n),i}_s \ud s} \\
		\nonumber
		&& + \sum_{i=1}^d sgn(c^{i,j}_n)N^{\{m-1\},i,j}_{|c^{i,j}_n|k^{m-1,(n),j}\delta\int_0^t Z^{m,(n),i}_r  \ud r},
		\qquad t\geq 0, \quad j \in D,
	\end{eqnarray*}
	where 
	\begin{equation*}
	X^{\{m-1\},(n),i,j}_r = X^{\{m-2\},(n),i,j}_{\int_{0}^{a_n\varepsilon} Z^{m-1,(n),i}_s \ud s + r} - X^{\{m-2\},(n),i,j}_{\int_{0}^{a_n\varepsilon} Z^{m-1,(n),i}_s \ud s}
	\end{equation*}
	and
	\begin{equation*}
	N^{\{m-1\},i,j}_r:=N^{\{m-2\},i,j}_{|c^{i,j}_n|k^{m-2,(n),j}\delta\int_{0}^{a_n\varepsilon} Z^{m-1,(n),i}_s \ud s + r} -N^{\{m-2\},i,j}_{|c^{i,j}_n|k^{m-2,(n),j}\delta\int_{0}^{a_n\varepsilon} Z^{m-1,(n),i}_s \ud s}. 
	\end{equation*}
	As before, we apply  Corollary 1 in \cite{CPGUB} to obtain that
	\begin{equation*}
	\left(\frac{a_n}{b^i_n}\VEC{Z}^{m,(n),i}_{a_nt}, t\geq 0, i\in D\right) \underset{n\rightarrow\infty}{\longrightarrow}
	\left(\VEC{Z}^{m,i}_t, t\geq 0, i \in D \right) 
	\end{equation*}
	in the Skorokhod space, where $\VEC{Z}^{m,i}$ is the unique solution of \eqref{eq: Zm}.

	Using the processes constructed above,  we define the process $\VEC{Z}^{\varepsilon,\delta,(n)}$ as
	\begin{equation*}
	Z^{\varepsilon,\delta,(n),j}_{t}:= Z^{\lceil s/a_n\varepsilon\rceil,(n), j}_{t/a_n\varepsilon - \lfloor t/a_n\varepsilon\rfloor}, \qquad t\geq 0 , j \in D,
	\end{equation*}
	which is the unique solution to \eqref{eq: DIMBP TC epsdelta}. 
	By the previous Skorokhod convergences, we can conclude that
	\begin{equation*}\label{eq: convz2}
	\left(\frac{a_n}{b^i_n}\VEC{Z}^{\varepsilon,\delta ,(n),i}_{a_nt},\ t\geq 0, i\in D\right) \underset{n\rightarrow\infty}{\longrightarrow}
	\left(\VEC{Z}^{\varepsilon, \delta,i}_t,\ t\geq 0, i \in D \right) 
	\end{equation*}
	in the Skorokhod space, where $\VEC{Z}^{\varepsilon,\delta}$ satisfies \eqref{eq: Zepsdelta}, and hence is the unique solution of \eqref{eq: Yepsdelta}.	
	\hfill \qedsymbol
	
	\subsection*{Proof of Lemma \ref{le: ZnZn,eps}}
	Define the process $\xi^{\varepsilon,\delta,(n),j}_{t}= Z^{(n),j}_{t} - Z^{\varepsilon,\delta,(n),j}_{t}$ and 
	$$\tau_M=\inf\{s\geq 0 : \tfrac{a_n}{b^i_n}Z_{a_ns}^{\varepsilon,\delta,(n),i}\geq M \text{ or } \tfrac{a_n}{b^i_n}Z^{(n),i}_{a_ns}\geq M \text{ for some } i \in \{1,\dots,d\}\}.$$
	We observe that $X^{(n),i,j}$ is a continuous time random walk for every $i,j$ and $n$. Therefore, it can be seen as a compound Poisson process, i.e.  
	\[
	X^{(n),i,j}_h=\int_{0}^{+\infty} \int_{-1}^{+\infty}\ell \indi{u \leq h} \mathcal{K}^{(n),i,j}(\ud \ell, \ud u),
	\] for some Poisson  random measure $\mathcal{K}^{(n),i,j}$ with intensity $\ud u \times \mu^{(n),i,j}(\ud\ell)$, where $\mu^{(n),i,j}$ is the jump distribution associated to $X^{(n),i,j}$. 
	Then, by  \eqref{eq: DIMBP TC epsdelta}  we have that
	\begin{eqnarray*}
		|\xi^{\varepsilon,\delta,(n),j}_{a_n(t \wedge \tau_M)}|
		&\leq&  \sum\limits_{i=1}^d \left|X^{(n),i,j}_{\int_0^{a_n(t \wedge \tau_M)}Z_s^{(n),i}\ud s} -  X^{(n),i,j}_{ \int_0^{a_n(t \wedge \tau_M)}Z_s^{\varepsilon,\delta,(n),i}\ud s}\right| \\ 
		&&+ \sum\limits_{i=1}^d \left|N^{i,j}_{|c_n^{i,j}|\int_0^{a_n(t \wedge \tau_M)}Z_s^{(n),i}Z_s^{(n),j}\ud s } - N^{i,j}_{|c_n^{i,j}|\sum\limits_{m=1}^{\lfloor t/\varepsilon\rfloor +1}\int_{a_n(m\varepsilon\wedge \tau_M)}^{a_n\left((m+1)\varepsilon \wedge t\wedge \tau_M\right)}Z_s^{\varepsilon,\delta,(n),i}\lfloor \tfrac{ Z_{ma_n\varepsilon}^{\varepsilon,\delta,(n),j}}{\delta}\rfloor \delta \ud s}\right|\\
		&\leq& \sum\limits_{i=1}^d \left|X^{(n),i,j}_{\int_0^{a_n(t \wedge \tau_M)}\left[Z_s^{(n),i} +\left(Z_s^{\varepsilon,\delta,(n),i} -Z_s^{(n),i}  \right)\right]\ud s } 
		-  X^{(n),i,j}_{ \int_0^{a_n(t \wedge \tau_M)}Z_s^{(n),i}\ud s} \right| 
		\indi{A^{(n),i}_{t}}\\
		&&+  \sum\limits_{i=1}^d \left|X^{(n),i,j}_{\int_0^{a_n(t \wedge \tau_M)}\left[Z_s^{\varepsilon,\delta,(n),i} - \left(Z_s^{\varepsilon,\delta,(n),i} -Z_s^{(n),i}  \right) \right]\ud s} -  X^{(n),i,j}_{ \int_0^{a_n(t \wedge \tau_M)}Z_s^{\varepsilon,\delta,(n),i}\ud s} \right| \indi{\left(A^{(n),i}_t\right)^c}\\ 
		&&+ \sum\limits_{i=1}^d \left|N^{i,j}_{|c_n^{i,j}|\int_0^{a_n(t \wedge \tau_M)}Z_s^{(n),i}Z_s^{(n),j}\ud s } - N^{i,j}_{|c_n^{i,j}|\sum\limits_{m=1}^{\lfloor t/\varepsilon\rfloor +1}\int_{a_n(m\varepsilon\wedge \tau_M)}^{a_n\left((m+1)\varepsilon \wedge t\wedge \tau_M\right)}Z_s^{\varepsilon,\delta,(n),i}\lfloor \tfrac{ Z_{ma_n\varepsilon}^{\varepsilon,\delta,(n),j}}{\delta}\rfloor \delta \ud s}\right|\\  
		&=& \sum\limits_{i=1}^d \left|X^{(n),i,j}_{\int_0^{a_n(t \wedge \tau_M)}Z_s^{(n),i}\ud s + \left|\int_0^{a_n(t \wedge \tau_M)}\xi^{\varepsilon,\delta, (n),i}_s \ud s\right| } 
		-  X^{(n),i,j}_{ \int_0^{a_n(t \wedge \tau_M)}Z_s^{(n),i}\ud s} \right|  \\
		&&+ \sum\limits_{i=1}^d \left|N^{i,j}_{|c_n^{i,j}|\int_0^{a_n(t \wedge \tau_M)}Z_s^{(n),i}Z_s^{(n),j}\ud s } - N^{i,j}_{|c_n^{i,j}|\sum\limits_{m=1}^{\lfloor t/\varepsilon\rfloor +1}\int_{a_n(m\varepsilon\wedge \tau_M)}^{a_n\left((m+1)\varepsilon \wedge t\wedge \tau_M\right)}Z_s^{\varepsilon,\delta,(n),i}\lfloor \tfrac{ Z_{ma_n\varepsilon}^{\varepsilon,\delta,(n),j}}{\delta}\rfloor \delta \ud s}\right|
	\end{eqnarray*}
	where 
	\begin{equation*}
	A^{(n),i}_{t}:=\left\{\int_0^{a_n(t \wedge \tau_M)}Z_s^{(n),i}\ud s  \leq \int_0^{a_n(t \wedge \tau_M)}Z_s^{\varepsilon,\delta,(n),i}\ud s\right\} 
	\end{equation*}
	with $\left(A^{(n),i}_{t}\right)^c$ its complement. By taking expectation we have
	\begin{equation*}
	\begin{split}
		\se\left[\left|\xi^{\varepsilon,\delta,(n),j}_{a_n(t \wedge \tau_M)} \right|\right]
		&\leq \sum\limits_{i=1}^d \se \left[ \left|X^{(n),i,j}_{\int_0^{a_n(t \wedge \tau_M)}Z_s^{(n),i}\ud s + \left|\int_0^{a_n(t \wedge \tau_M)}\xi^{\varepsilon,\delta, (n),i}_s \ud s\right| } 
		-  X^{(n),i,j}_{ \int_0^{a_n(t \wedge \tau_M)}Z_s^{(n),i}\ud s} \right| 
		\right]\\
		&\hspace{-.7cm}+ \sum\limits_{i=1}^d \se \left[ N^{i,j}_{\left|c_n^{i,j}\sum\limits_{m=1}^{\lfloor t/\varepsilon\rfloor +1}\int_{a_n(m\varepsilon\wedge \tau_M)}^{a_n\left((m+1)\varepsilon \wedge t\wedge \tau_M\right)} \left( Z_s^{(n),i}Z_s^{(n),j} - 
			Z_s^{\varepsilon,\delta,(n),i}\lfloor \tfrac{ Z_{ma_n\varepsilon}^{\varepsilon,\delta,(n),j}}{\delta}\rfloor \delta\right)\ud s\right|}\right]\\
		&\hspace{-.8cm}\leq  \sum\limits_{i=1}^d \se \left[\left|X^{(n),i,j}_{\left|\int_0^{t}
			a_n \xi^{\varepsilon,\delta,(n),i}_{a_n(s \wedge \tau_M)} \ud s\right|}\right|\right]\\
		&\hspace{-.7cm}+ \sum\limits_{i=1}^d |c_n^{i,j}|\sum\limits_{m=1}^{\lfloor t/\varepsilon\rfloor +1}\int_{m\varepsilon}^{(m+1)\varepsilon \wedge t}\se \left[ a_n\left|  Z_{a_n(s \wedge \tau_M)}^{(n),i} Z_{a_n(s \wedge \tau_M)}^{(n),j} - 
		Z_{a_n(s \wedge \tau_M)}^{\varepsilon,\delta,(n),i}\lfloor \tfrac{ Z_{a_n(m\varepsilon \wedge \tau_M)}^{\varepsilon,\delta,(n),j}}{\delta}\rfloor \delta \right|\right]\ud s 	 
		\end{split}
	\end{equation*}	
	Let
	\begin{equation*}
	S^{(n),i,j}_h=\int_{0}^{+\infty} \int_{-1}^{1}\ell \indi{u \leq h} \mathcal{K}^{(n),i,j}(\ud \ell, \ud u)\qquad \mbox{and} \qquad R^{(n),i,j}_h=\int_{0}^{+\infty} \int_{1}^{+\infty}\ell \indi{u \leq h} \mathcal{K}^{(n),i,j}(\ud \ell, \ud u),
	\end{equation*}
	Then,
	\begin{eqnarray*}
		\se\left[\left|\xi^{\varepsilon,\delta,(n),j}_{a_n(t \wedge \tau_M)} \right|\right]
		&\leq&  \sum\limits_{i=1}^d \se \left[\left|S^{(n),i,j}_{\left|\int_0^{t}
			a_n\xi^{\varepsilon,\delta,(n),i}_{a_n(s \wedge \tau_M)} \ud s\right|}\right|\right] + \sum\limits_{i=1}^d \se \left[R^{(n),i,j}_{\left|\int_0^{t}
			a_n\xi^{\varepsilon,\delta,(n),i}_{a_n(s \wedge \tau_M)} \ud s\right|}\right]  
		\\
		&+& \sum\limits_{i=1}^d  |c_n^{i,j}| \int_{0}^{t}\se\left[ a_nZ^{(n),j}_{a_n(s\wedge \tau_M)}\left|\xi_{a_n(s\wedge \tau_M)}^{\varepsilon,\delta,(n),i}\right| + a_nZ^{\varepsilon,\delta,(n),i}_{a_n(s\wedge \tau_M)}\left|\xi_{a_n(s\wedge \tau_M)}^{\varepsilon,\delta,(n),j}\right| \right] \ud s\\
		&+& \sum\limits_{i=1}^d |c_n^{i,j}|\sum\limits_{m=1}^{\lfloor t/\varepsilon\rfloor+1}\int_{m\varepsilon}^{(m+1)\varepsilon \wedge t} \se\left[a_nZ^{\varepsilon,\delta,(n),i}_{a_n(s\wedge \tau_M)}\left|Z_{a_n(s\wedge \tau_M)}^{\varepsilon,\delta,(n),j} - \lfloor \tfrac{ Z_{a_n(m\varepsilon\wedge \tau_M)}^{\varepsilon,\delta,(n),j}}{\delta}\rfloor \delta\right|\right]\ud s,	 
	\end{eqnarray*}	
	so taking in account the scale constant we have
	\begin{eqnarray*}
		\se\left[\tfrac{a_n}{b^j_n}\left|\xi^{\varepsilon,\delta,(n),j}_{a_n(t \wedge \tau_M)} \right|\right]
		&\leq& \sum\limits_{i=1}^d \tfrac{a_nb^i_n}{b^j_n}\left[\mu^{(n),i,j}(\{-1,1\})+ \int_1^{+\infty}\ell \mu^{(n),i,j}(\ud \ell)\right]\int_0^{t}\se \left[
		\left|\tfrac{a_n}{b^i_n}\xi^{\varepsilon,\delta,(n),i}_{a_n(s \wedge \tau_M)}\right| \right]\ud s\\
		&+& \sum\limits_{i=1}^d  |c_n^{i,j}|b_n^i \int_{0}^{t}\se\left[ \tfrac{a_n}{b_n^j}Z^{(n),j}_{a_n(s\wedge \tau_M)}\tfrac{a_n}{b_n^i}\left|\xi_{a_n(s\wedge \tau_M)}^{\varepsilon,\delta,(n),i}\right| + \tfrac{a_n}{b_n^i}Z^{\varepsilon,\delta,(n),i}_{a_n(s\wedge \tau_M)}\tfrac{a_n}{b_n^j}\left|\xi_{a_n(s\wedge \tau_M)}^{\varepsilon,\delta,(n),j}\right| \right] \ud s\\
		&+& \sum\limits_{i=1}^d |c_n^{i,j}|b_n^i\sum\limits_{m=1}^{\lfloor t/\varepsilon\rfloor+1}\int_{m\varepsilon}^{(m+1)\varepsilon \wedge t} \se\left[\tfrac{a_n}{b_n^i}Z^{\varepsilon,\delta,(n),i}_{a_n(s\wedge \tau_M)}\tfrac{a_n}{b_n^j}\left|Z_{a_n(s\wedge \tau_M)}^{\varepsilon,\delta,(n),j} - \lfloor \tfrac{ Z_{a_n(m\varepsilon\wedge \tau_M)}^{\varepsilon,\delta,(n),j}}{\delta}\rfloor \delta\right|\right]\ud s\\
		 & \leq & \sum\limits_{i=1}^d \tfrac{a_nb^i_n}{b^j_n}\left[\mu^{(n),i,j}(\{-1,1\})+ \int_1^{+\infty}\ell \mu^{(n),i,j}(\ud \ell)\right]\int_0^{t}\se \left[
		\left|\tfrac{a_n}{b^i_n}\xi^{\varepsilon,\delta,(n),i}_{a_n(s \wedge \tau_M)}\right| \right]\ud s\\
		&+& \sum\limits_{i=1}^d  M|c_n^{i,j}|b_n^i \int_{0}^{t}\se\left[ \tfrac{a_n}{b_n^i}\left|\xi_{a_n(s\wedge \tau_M)}^{\varepsilon,\delta,(n),i}\right| + \tfrac{a_n}{b_n^j}\left|\xi_{a_n(s\wedge \tau_M)}^{\varepsilon,\delta,(n),j}\right| \right] \ud s\\
		&+& \sum\limits_{i=1}^d M|c_n^{i,j}|b_n^i\sum\limits_{m=1}^{\lfloor t/\varepsilon\rfloor+1}\int_{m\varepsilon}^{(m+1)\varepsilon \wedge t} \se\left[\tfrac{a_n}{b_n^j}\left|Z_{a_n(s\wedge \tau_M)}^{\varepsilon,\delta,(n),j} - \lfloor \tfrac{ Z_{a_n(m\varepsilon\wedge \tau_M)}^{\varepsilon,\delta,(n),j}}{\delta}\rfloor \delta\right|\right]\ud s.	 
	\end{eqnarray*}
	Note that by the scaling limit hypothesis \eqref{eq: CL},   we have 
	\begin{eqnarray*}
		\tfrac{a_nb^i_n}{b^j_n}	\mu^{(n),i,j}(\{-1,1\}) &+& 	\tfrac{a_nb^i_n}{b^j_n}\int_1^{+\infty}\ell \mu^{(n),i,j}(\ud \ell)\\
		&=&
		2\tfrac{a_n}{b^j_n}\left[\int_0^{b^i_n}\int_{-1}^{+\infty}\indi{\{\ell=-1\}}\mu^{(n),i,j}(\ud \ell)\ud s\right]
		+ \tfrac{a_n}{b^j_n}\int_0^{b^i_n}\int_{-1}^{+\infty}\ell \ud s\mu(\ud \ell )\\ 
		&=&2\se\left[\tfrac{a_n}{b^j_n}\int_0^{b^i_n}\int_{-1}^{+\infty}\indi{\{\ell=-1\}}\mathcal{K}(\ud \ell,\ud s)\right] +\se\left[\tfrac{a_n}{b^j_n}X^{(n),i,j}_{b^i_n}\right]
		\underset{n \rightarrow \infty}{\longrightarrow}  \se\left[X^{i,j}_{1}\right],
	\end{eqnarray*}
	where we use that
	$\se\left[\tfrac{a_n}{b^j_n}\int_0^{b^i_n}\int_{-1}^{+\infty}\indi{\{\ell=-1\}}\mathcal{K}(\ud \ell,\ud s)\right]$  is the (reescaled) number of negative jumps of the process $X^{(n),i,j}$ in the interval $[0,b_n^i]$, that goes to zero in the limit. Therefore, for $n$ big enough  
	\begin{equation*}
	\max\limits_{i,j \in \{1,\dots,d\}}	\left\{\tfrac{a_nb^i_n}{b^j_n}	\mu^{(n),i,j}(\{-1,1\}) + 	\tfrac{a_nb^i_n}{b^j_n}\int_1^{+\infty}\ell \mu^{(n),i,j}(\ud \ell)\right\} \leq 1 + \max\limits_{i,j \in \{1,\dots,d\}} \left|\se\left[X^{i,j}_{1}\right]\right|=: I ,
	\end{equation*}
	and if we take $C_n =\max\limits_{i,j \in \{1,\dots,d\}}|c^{i,j}_nb^i_n|$, we have 
	\begin{eqnarray*}
		\sum\limits_{j=1}^d	\se\left[\left|\tfrac{a_n}{b^j_n}\xi^{\varepsilon,\delta,(n),j}_{a_n(t \wedge \tau_M)} \right|\right]
		&\leq&  d\left(I + 2MC_n\right) \int_0^t \sum\limits_{i=1}^d \se \left[
		\left|\tfrac{a_n}{b^i_n}\xi^{\varepsilon,\delta,(n),i}_{a_n(s \wedge \tau_M)}\right| \right]\ud s \\
		&&+ dMC_n \sum\limits_{j=1}^d \sum\limits_{m=1}^{\lfloor t/\varepsilon\rfloor+1}\int_{m\varepsilon}^{(m+1)\varepsilon \wedge t} \se\left[\tfrac{a_n}{b_n^j}\left|Z_{a_n(s\wedge \tau_M)}^{\varepsilon,\delta,(n),j} - \lfloor \tfrac{ Z_{a_n(m\varepsilon\wedge \tau_M)}^{\varepsilon,\delta,(n),j}}{\delta}\rfloor \delta\right|\right]\ud s.	 
	\end{eqnarray*}
	
	Therefore, by Gronwall's Lemma,
	\begin{eqnarray*}
		&\sum\limits_{j=1}^d&\se\left[\left|\tfrac{a_n}{b^j_n}\xi^{\varepsilon,\delta,(n),j}_{a_n(t \wedge \tau_M)}\right|\right]\\
		&\leq&  e^{td\left[I + 2MC_n\right]}dMC_n\sum\limits_{j=1}^d \sum\limits_{m=1}^{\lfloor t/\varepsilon\rfloor+1} \int_{m\varepsilon}^{(m+1)\varepsilon \wedge t}\se\left[  \tfrac{a_n}{b^j_n}\left|Z_{a_n(s\wedge \tau_M)}^{\varepsilon,\delta,(n),j} - \lfloor \tfrac{ Z_{a_n(m\varepsilon\wedge \tau_M)}^{\varepsilon,\delta,(n),j}}{\delta}\rfloor \delta\right| \right]\ud s.
	\end{eqnarray*}
	
	By taking the upper limit as $n$ goes to infinity and using Lemma \ref{le: SL epsdelta}, we have that
	\begin{equation*}
	\limsup_{n}\sum\limits_{j=1}^d\se\left[\left|\tfrac{a_n}{b^j_n}\xi^{\varepsilon,\delta,(n),j}_{a_n(t \wedge \tau_M)}\right|\right] 
	\leq  e^{td\left[I + 2MC\right]}dMC\sum\limits_{j=1}^d \sum\limits_{m=1}^{\lfloor t/\varepsilon\rfloor+1} \int_{m\varepsilon}^{(m+1)\varepsilon \wedge t}\se\left[  \left|Z_{s}^{\varepsilon,\delta,j} - \lfloor \tfrac{ Z_{m\varepsilon}^{\varepsilon,\delta,j}}{\delta}\rfloor \delta\right| \right]\ud s
	\end{equation*}
	where $C=\max\limits_{i,j \in \{1,\dots,d\}}|c^{i,j}|$.  Finally, by the c\'adl\'ag properties of the process $\VEC{Z^{\varepsilon,\delta}}$, we have that
	\begin{equation*}
	\lim\limits_{\varepsilon,\delta \rightarrow 0}\lim\limits_{n\rightarrow \infty}\sum\limits_{j=1}^d\se\left[\left|\tfrac{a_n}{b^j_n}\xi^{\varepsilon,\delta,(n),j}_{a_n(t \wedge \tau_M)}\right|\right]=0
	\qquad \forall t\geq 0 ,
	\end{equation*}
	and we deduce that $((Z^{(n),i} - Z^{\varepsilon,\delta,(n),i}), i \in \mathcal{D}) \rightarrow 0$ in the sense of finite dimensional distribution. Using the uniqueness of the processes, we deduce the almost sure convergence in Skorokhod space.
	
	\hfill \qedsymbol

\section*{Acknowledgements}
Both authors would like to thank Ger\'onimo Uribe Bravo for collaborate on an early draft of this paper. 
SP's research is supported by PAPIIT IA103220.

\bigskip

\bibliographystyle{plain}

\end{document}